\documentclass{article}

\usepackage{amsmath, amsfonts, amsthm, amssymb, amscd}
\usepackage[mathcal]{euscript} 
\usepackage{mathrsfs}  
\newtheorem{theorem}{Theorem}[section]

\newtheorem{definition}[theorem]{Definition} 
\newtheorem{remark}[theorem]{Remark}
\numberwithin{equation}{section} 

\usepackage[all]{xy}
\usepackage{amscd}
\usepackage{amsthm}
\usepackage{graphicx}
\usepackage{color}
\usepackage{hyperref}

\def\f{\lambda}
\def\g{\mu}

\def\Rp{{R}^+}
\def\bft{\mathbf{t}}

\def\BI{Moore(B)}

\begin{document}

\title{``Parallel'' transport - revisited}
\author{jim stasheff}

\maketitle

\begin{abstract}
Parallel transport in a fibre bundle with respect to smooth paths in the base space$B$ have recently been extended to  representations of the smooth singular simplicial set $Sing_{smooth}(B)$.
  Inspired by these extensions,
 I revisit the development of a notion of `parallel' transport in the topological setting of fibrations with the homotopy lifting property and then extend it to representations of $Sing(B)$ on  such fibrations. Closely related is the notion of (strong or $\infty$) homotopy action,
which has variants under a variety of names. 
\end{abstract}


\section{Introduction/History}
In classical differential geometry
(``a language the muse did not sing at my cradle'' - see below), \emph{parallel transport} is defined in the context of a \emph{connection} on a smooth bundle  $p:E\to B$.
The latter can mean a covariant derivative operator, a differential 1-form or a set of horizontal subspaces in the tangent bundle $Tp:TE\to TB.$ The corresponding \emph{parallel transport} 
$\tau: E\times_B  B^I \to E$ is constructed by lifting a path in $B$ to a \emph{unique!} path in $E$ with specified starting point. 
The \emph{holonomy} is given by the evaluation of $\tau$ on $\Omega B$, the space of based loops in $B$. The \emph{holonomy group} is the image as a subgroup of the structure group of the bundle. That it is a group follows from the uniqueness of the lifting. It is well defined up to conjugation depending on the choice of base point.

If $p:E\to B$ is only a fibration of topological spaces, the situation is different: we still can lift paths but not uniquely.
However, the lifts are related by homotopies, so that the  right action of 
$\Omega B$ on the fibre $F$ is well behaved up to homotopy, i.e.

\begin{displaymath}
\begin{CD}
   F\times  \Omega B \times    \Omega B @>m \times 1>> F \times \Omega B \\
    @V{1\times a}VV              @VV{a}V \\
F \times \Omega B    @> a>>    F
\end{CD}
\end{displaymath}

\noindent is homotopy commutative, where $m$ is multiplication in $\Omega B$ and $a$ is the action. In fact, there is a whole sequence of coherent higher homotopies, which we review below.
(Because $\lambda\mu$ denotes travelling along $\f$ first and then along $\g$, the action will be written as right action: $(f,\f)\mapsto f\f.$)

Perhaps the oldest treatment in algebraic topology
(I learned it as a grad student from Hilton's Introduction to Homtopy Theory \cite{hilton:intro}  - the earliest textbook on the topic)
is to  consider the long exact sequence, where $F$ is the fibre over a chosen base point in $B$, 
$$\cdots\to \pi_n(F)\to \pi_n(E)\to \pi_n(B)\to \pi_{n-1}(F)\to \cdots$$
ending with $$\cdots \to \pi_1(B)\to \pi_0(F)\to \pi_0(E)\to \pi_0(B).$$
Of course, exactness is very weak at the end since the last three are in general only sets, but
exactness at $ \pi_0(F)$ is in terms of the action of $ \pi_1(B)$ on $ \pi_0(F).$  This passage to homotopy classes obscures the `action' of $\Omega B$ on $F$.  Initially, this was referred to as a \emph{homotopy action},
meaning only that $\lambda(\mu f)$ was homotopic to $(\lambda\mu)f$  for $f\in F$ and $\lambda, \mu \in \Omega B.$

In  those days, at least at Princeton, there was no differential geometry until Milnor gave an undergrad course my final year there.  Notice that our book `Characteristic Classes' \cite{MS:cc} considers differential forms only in Appendix C, added much later. I think this lack of de Rham theory in my graduate education was the result of Serre's thesis which triumphed over characteristic 0, cf. Colloque de Topologie, Bruxelles (1950). Only recently have I learned of even greater significance of the metric geometry for the boundary map $ \pi_n(B) \to \pi_{n-1}(F)$ \cite{duran:exotic}.

Similarly,  I learned only much lateer of the notion of \emph{thin homotopy} which quotients $\Omega B$ to a group without losing so much information. Just recently, Johannes Huebschman led me to a paper of Kobayashi (from 1954!) where he is already using what is now called thin homotopy in terms of parallel transport and holonomy for smooth bundles with connection.

Back in 1966,  in the Mexican Math Bulletin \cite{jds:mex}, a journal not readily available\cite{jds:mex} (although summarized in \cite{jds:madison,jds:bingh}),  I showed that in the topological setting of fibrations 
(satisfying  the homotopy lifting property), there was a notion of `parallel'  transport not dependent on having a connection. This meant  not only the above homotopy action, but in fact an sh (or 
$A_\infty$)-action, which is to say the adjoint $\Omega B\to End(F)$ was an 
$A_\infty$-map.
 (Not long after, Nowlan \cite{nowlan:action} studied more general $A_\infty$-actions of fibres on total spaces.) Note that the homotopy lifting property does \emph{not} imply unique lifting, so even the composition of paths is \emph{not} necessarily respected when lifting.

For my purposes, it was sufficient to consider transport along based loops in the base,  though the arguments  allow for transport along any path in the base.

Inspired by Block-Smith \cite{block-smith} and Igusa  \cite{igusa:twisting}, Abad and Schaetz \cite{abad-schaetz}  showed that `flat' smooth parallel transport (e.g. as in  \cite{igusa:twisting,block-smith}) can be derived from the $A_\infty$ version of de Rham's theorem due to Gugenheim \cite{gugenheim:deRham}; this allows them to extend parallel transport to 
an $A_\infty$-functor.\footnote{The term `flat'  is used to denote  direct generalization of the usual notion of parallel transport with zero curvature. More  generally, smooth (not necessarily flat) parallel transport may be described in terms of integration of $\infty$-connections data, as in \cite{fss}.}

They work with what is known as the  $\infty$-groupoid $\Pi_{\infty}(B)$ of a space $B$ and 
its \emph{representations up to homotopy} on a fibre bundle $E\to B$. The   $\infty$-groupoid $\Pi_{\infty}(B)$ can be represented as a simplicial set: the singular complex $Sing(B).$ This led me to revisit my earlier work on `parallel' transport and extend   work to the fibrations setting, i.e. without any smoothness or connection operator or differential form or any infinitesimal input.

 What I am  after here is an analog of parallel transport along paths extended to parallel transport over (maps of) simplices, not  just 1-simplices.


\begin{remark} Just as there is only linguistic difference between an action of a group on an object and a representation of that group, the same is true in the homotopy setting, but sometimes a homotopy action means only the existence of a single homotopy and some times the usual coherent collection of higher homotopies;
when in doubt, it is best to prepend $A_\infty$ or \emph{strong homotopy}.
\end{remark}
\section{The `classical' topological case}
We first recall what are rightly known as Moore paths \cite{jcm:paths} on a topological space $X$.
\begin{definition}
Let  $R^+=[0, \infty)$ be the nonnegative real line. For a space $X,$
let $\emph{Moore(X)}$ be the subspace of \emph{Moore paths}  $\subset X^{\Rp}\times \Rp$ of pairs $(f,r)$
such that $f$ is constant on $[r, \infty)$. There are two
maps

\begin{itemize}
\item
$\partial^-, \partial^+: Moore(X) \to X,$
\item
$\partial^-(f,r)=f(0),$
\item
$\partial^+(f,r)=f(r).$
\end{itemize}

Recall composition $\circ$ of Moore paths in $Moore(X)$ is given by 
sending pairs
$(\lambda,r),(\g,s)\in Moore(X)$ such that $\f(r)=\g(0)$ to  $\f\g\in Moore(X) $ which  is constant on $[r+s,
\infty)$, $\f\g|[0,r]=\f|[0,r]$ and $\f\g(t)=\g(t-r)$ for $t \geq r$. 
An identity function $\epsilon: X \to Moore(X)$ is given by
$\epsilon(x)=(\hat{x},0)$ where $\hat{x}$ is the constant map on
$\Rp$ with value $x$.
\end{definition}
Composition is continuous and 
 gives, as is well known, a category/$\infty$-groupoid structure on 
$Moore(X)$. 

 If we had used the `ancient' Poincar\'e paths $I\to X$, we would have had to work with an $A_\infty$-structure on $X^I.$ Indeed, it was working with that standard parameterization which led to $A_\infty$-structures
\cite{sugawara:g, jds:hahI}.

For a category $C$, we denote by $C_{(n)}$ the set of n-tuples of composable morphisms.  In partcular, we will be concerned with $\BI_{(n)}.$ We will write $\mathbf{t}$ for $(t_1,\cdots,t_n)$ and $\hat t_i$ for $(t_1,\cdots,t_{i-1}, t_{i+1},\cdots,t_n)$,    
 Back in 1988 \cite{jds:shrep}, I referred to \emph{strong homotopy representations}, but today I will use the \emph{representation up to homotopy} terminology, having in mind the generalization that comes next.

\begin{definition}
\label{Thetas}
A \emph{representation up to homotopy} of  $\BI$ on a fibration $E\to B$ is 
an $A_\infty$-morphism (or shm-morphism \cite{sugawara:hc})  from $\BI$ to $End_B(E)$; that is, a collection of maps
 $$\theta_n: I^{n-1}\times E  \times_B \BI_{(n)}  \to E$$
 (where $E  \times_B \BI_{(n)} $ consists of $n+1$-tuples  $(e,\lambda_1,\dots, \lambda_n)$ where the
 $\lambda_i$  are composable paths, constant on $[r_i,\infty)$, and $p(e) = \lambda_1(0)$)
 such that
$$p(\theta_n(\bft,e,\lambda_1,\dots, \lambda_n)) = \lambda_n (r_n),$$
$$\theta_n(\mathbf{t},--,\lambda_1,\dots, \lambda_n)$$
is a fibre homotopy equivalence and
 satisfies the usual/standard relations:
 \begin{itemize}
\item $ \; 
\displaystyle{\theta_n(t_1,\cdots,t_i=0,\cdots, t_{n-1},e,\lambda_1,\dots, \lambda_n)) 
= \theta_{n-1}(\hat t_i,e,\cdots,\lambda_{i}\lambda_{i+1},\cdots))}$

\item $ \;
\displaystyle{\theta_n(t_1,\cdots,t_i=1,\cdots, t_{n-1},e,\lambda_1,\dots, \lambda_n)) = } \\
\hspace*{\stretch{1}}
\displaystyle{ = \theta_i(\cdots,t_{i-1}, \theta_{n-i}(t_{i+1},\cdots,t_{n-1}, \lambda_i,\cdots, \lambda_n,e),\lambda_1,\dots, \lambda_{i-1},).} 
\hspace*{\stretch{0}}$
 \end{itemize}
 \end{definition}
 \begin{remark} \emph{That the parameterization is by cubes, as for Sugawara's strongly homotopy multiplicative maps rather than more general polytopes, reflects the fact that $Moore(X)$ and $ End_B(E)$ are strictly associative. Strictly speaking, referring to $\BI \to End_B(E)$ as an $A_\infty$-map raises issues about a topology on  $End_B(E);$ the adjoint formulas above avoid this difficulty.
 \vskip1ex
 Even more difficult, at least for exposition and detailed proofs, would be the use of `Poincar\'e'  loops requiring $\theta_n$ to use $K_{n+1}$ instead of $ I^{n-1}.$}
   \end{remark}

 Since our construction uses in a crucial way the homotopy lifting property, we first construct maps
$$\Theta_n: I^n \times E  \times_B \BI_{(n)}   \to E$$
such that the desired $\theta_n$ are then recovered at $t_1=1.$

 \begin{remark} \emph{The special case} $\Theta_1: I \times E  \times_B \BI_{(n)}  \to E$ \emph{corresponds to a path lifting function. In a little known paper \cite{d&e}, Dyer and Eilenberg call such a $\Theta_1$ an action of the path space $\BI_{(1)}$ on $E$ and point out that its existence is equivalent to $p$ having the homotopy lifting property.}
\end{remark}

The idea is that
 if $\Theta_j$ has been defined satisfying these relations for all $j<n$, the $\Theta_{n-1}$ will fit together to define $\Theta_n$ on all faces of the cube except for the face where $t_1=1.$ In analogy with the horns of simplicial theory, we will talk about filling an \emph{open  box}, meaning the boundary of the cube minus the open face, called a \emph{lid}, where $t_i=1$ (compare horn-filling in the simplicial setting).  Use the homotopy lifting property to  `fill in the box'  in $E$ after filling in the trivial image box in $B.$ That box in $B$ is  a trivial box since it is just  the composite path $\lambda_1\cdots \lambda_n.$

It might help to consider the cases $n=1,2.$
Consider $(\lambda,r)\in Moore(B).$ Lift $\lambda$ to a path $(\bar\lambda,r)$ starting at $e\in E.$
Define $\Theta_1: I \times E\to Moore (B) \to E$ by

$$\Theta_1(t,e,(\lambda,r)) = (\bar\lambda,r)(tr)\in E$$
and $\theta_1(e,(\lambda,r)) = \Theta_1(1,e,(\lambda,r)) =: e(\lambda,r).$ 

Now lift $(\mu,s)$ to  a path $(\bar\mu,s)$ starting at $e(\lambda,r)\in E$
and lift $ (\lambda,r)(\mu,s)$ to a path $(\overline{\lambda\mu},r+s)$ starting at $e.$
These lifts fit together to define a map to $E$, which will be the restriction of the desired map   on the open 2-dimensional box of the desired map $\Theta_2.$ This open box has an  image in $B$ which can trivially be filled in. Regarding the filling as a homotopy, the map to $E$ on the open 2-dimensional box can be filled in by lifting that homotopy.

\begin{theorem}
\label{parallel}
(cf. Theorem A in \cite{jds:mex}) For any fibration $p:E\to B$, there is an $A_\infty$-action $\{\theta_n\}$ of Moore(B) on $E$ such that $\theta_1$ is a fibre homotopy equivalence. This action is unique up to homotopy
in the  $A_\infty$-sense.

\end{theorem}

In Theorem B in \cite{jds:mex}, I proved further:

\begin{theorem}
Given an $A_\infty$-action $\{\theta_n\}$ of the Moore loops $\Omega B$ on a space $F$, there is a fibre space $p_\theta: E_\theta\to B$ such that, up to homotopy,  the $A_\infty$-action $\{\theta_n\}$ can be recovered by the above procedure. If the $A_\infty$-action $\{\theta_n\}$ was originally obtained by the above procedure from a fibre space $p:E\to B$, then $p_\theta$ is fibre homotopy equivalent to $p$.
\end{theorem}

This construction gave rise to the slightly more general (re)construction below.  It can also be generalized to give an $\infty$-version of the Borel construction/homotopy quotient: $G\to X \to X_G=X//G$ for  an $A_\infty$-action \cite{iwase-mimura:hha}.
 
\section{Upping the ante to {\it Sing}}
 Inspired by Block-Smith \cite{block-smith} and Igusa  \cite{igusa:twisting}, Abad and Schaetz \cite{abad-schaetz}   look not at just composable paths, but rather  look at the singular complex $Sing(B)$. For a singular $k$-simplex $\sigma:\Delta^k\to B$, there are several $k$-tuples  of composable paths  from vertex 0 to vertex $k$ by restriction to edges, in fact, $k!$ such.  Given $\sigma$, we denote by $F_i$ the fibre over vertex $i\in \sigma$.
 
  Following e.g.\! Abad-Schaetz \cite{abad-schaetz} (based on Abad's thesis and his earlier work with Crainic), we make the following definition of a \emph{representation up to homotopy}, where we take  a singular $k$-simplex $\sigma$ to be (the image of )\newline $<0,1,\cdots,k>$ with the $p$-th face $\partial_p\sigma$ being $<0,\cdots,p-1,p+1,\cdots,k>.$ However, we keep much of the notation above rather than switch to theirs. 
  
  \begin{definition}
 A \emph{representation up to homotopy} of $Sing(B)$ on a fibration $E\to B$ is a collection of maps
$\{ \theta_{k}\}_{k\geq 0}$ which assign to any $k$-simplex $\sigma:\Delta^k\to B$ a map $ \theta_{k}(\sigma): I^{k-1}\times F_0 \to F_k$ satisfying, for any $e\in F_0$, the relations:

$\theta_0$ is the identity on  $F_0$

\noindent For any $(t_1,\cdots,t_{k-1})$, 

$\theta_k(\sigma)(t_1,\cdots,t_{k-1}, -): F_0 \to F_k$ is a homotopy equivalence.

\noindent For any $1\leq p\leq k-1$ and $e\in F_0$, 

\noindent

$$ \theta_{k}(\sigma)(\cdots ,  t_p=0, \cdots, e) =  \theta_{k-1}(\partial_p\sigma)(\cdots,\hat t_p,\cdots,, e)$$

$$ \theta_{k}(\sigma)(\cdots ,  t_p=1, \cdots, e) = $$ 
$$
 \theta_{p}(<0,\cdots,p>)(t_1,\cdots, t_{p-1},  \theta_{q}(<p,\cdots,k>)(t_{p+1},\cdots,t_k,e)).$$

\end{definition}


\begin{remark} \emph{In definition \ref{parallel}, we worked with Moore paths so that the $A_\infty$- map was between
strictly associative spaces. Here instead the compatible 1-simplices compose just as e.g. a pair of 1-simplices and are related to a single 1-simplex only by an intervening 2-simplex. Associativity is trivial; the subtlety is in handling the 2-simplices and higher ones for multiple compositions. The idea of constructing a representation up to homotopy is very much like that of Theorem 1, the major difference being that instead of comparing two different liftings of the composed paths which are necessarily homotopic, we are comparing a  lifting e.g. of a path from 0 to 1 to 2 with a lifting of a path from 0 to 2 IF there is a singular 2-simplex $<012>$.  However, note that $<02>$ plays the role of $\lambda_1\lambda_2$ of Moore paths in the above formulas.}
\end{remark}

\begin{theorem}For any fibration $p:E\to B$, there is a representation up to homotopy of $Sing(B)$ on $E.$ 
\end{theorem}

\begin{remark} \emph{The fact that the representation up to homotopy is by fibre homotopy equivalence (as for the action of Moore(B)) is justified by the following:
Since a simplex $\sigma$ is contractible, the pullback $\sigma^*E$ is fibre homotopy trivial over $\sigma$.
Choose the requisite lifts in $\sigma^*E$ using a trivialization corresponding to the homotopy we want to lift and then map back into $E$.} 
\end{remark}
\begin{remark} \emph
{In contrast to the smooth bundle case where a connection provides unique path lifting, the fibration case is considerably more subtle since path and homotopy lifting is far from unique.
}
\end{remark}
The proof is in essence the same as that for Theorem \ref{parallel}.  The desired $\theta_n$ will appear as the missing lid on an open box (defined inductively) which is filled in by homotopy liftings $\Theta_n$ of a \emph{coherent} set of maps
$$p_n:I^{n}\to \Delta^n,$$
  where  $\Delta^n$ is the set
$$\{(t_1,\cdots, t_n) |  0\leq t_1 \leq t_2...\leq 1\},$$ are given in terms of 
 iterated convex linear functions. 
 
\begin{figure}[htbp]
\begin{center}
\includegraphics[width=.5\textwidth]{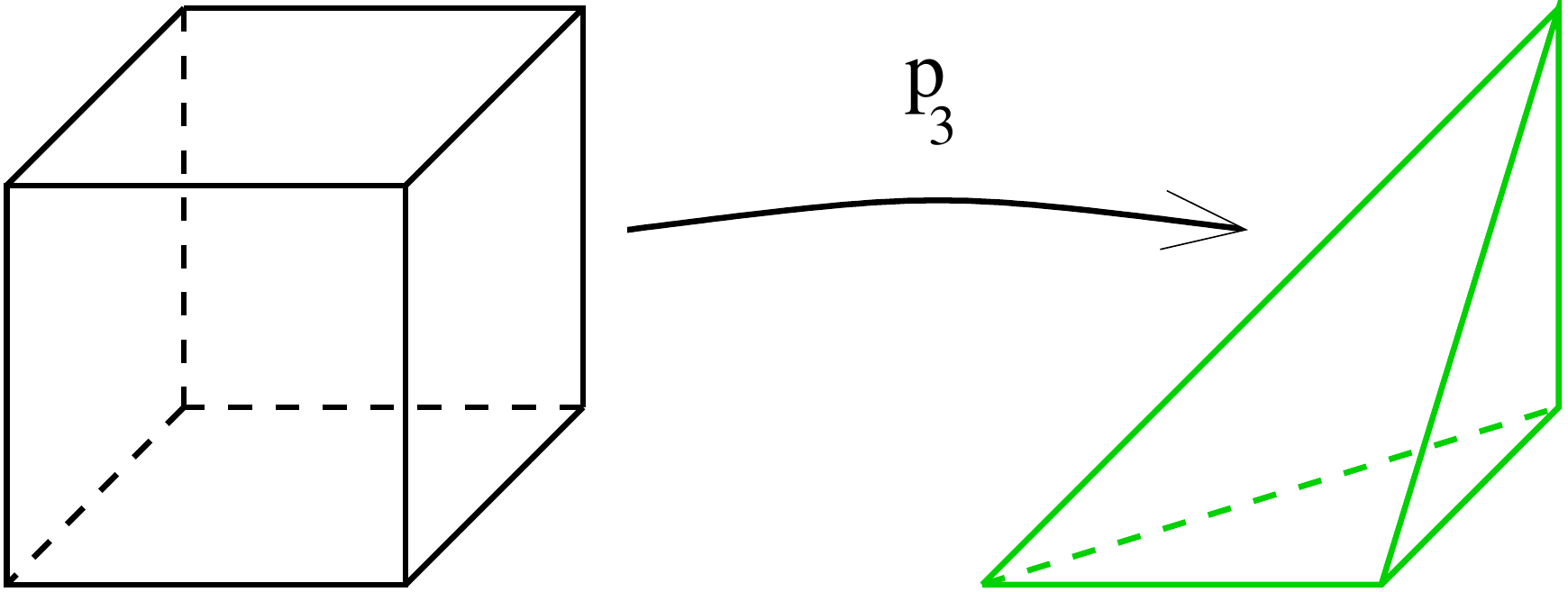}		
\caption{$p_3:I^3\to \Delta^3$}
\label{cube2simplex}
\end{center}
\end{figure}

 For $n=1$, define $p_1: t\mapsto t_1.$
 For $n=2,3$, write $t_1=t, \ t_2 =s,\  t_3=r.$
For $n=2$, define $c:=p_2:(t,s)\mapsto (t\cdot 1+ (1-t)s,s)$ and then 
$$p_3:(t,s,r)\mapsto (c(c(t,s),r), c(s,r),r) = (c(t\cdot 1+ (1-t)s),r),c(s,r),r).$$
See Figure \ref{cube2simplex}.

These have probably been written elsewhere; if you find them, let me know.

By coherent, I mean respecting the facial structure of the cubes and simplices.
Closely related are coherent maps
$$\gamma_n:I^{n-1}\to P\Delta^n$$
where $P$ denotes the set of paths, i.e. $P\Delta^n= Map(I,\Delta^n).$ 
Such maps were first produced by Adams \cite{adams:cobar} in the topological context by induction using the contractability of $\Delta^n$.  Later specific formulas were introduced by Chen \cite{chen:iterated} and, most recently, equivalently but more  transparently, by Igusa \cite{igusa:twisting}. 

More precisely:

\quad\quad\quad\quad\quad\quad\quad $\gamma_1(0)$ is the trivial path, constant at $0$

\quad\quad\quad\quad\quad\quad\quad and $\gamma_2:I \to \Delta^1$ is the `identity'. 

\noindent For any $1\leq p\leq k-1$,

$$\gamma_{k}(\cdots ,  t_p=0, \cdots) = \gamma_{k-1}(\cdots,\hat t_p,\cdots)$$
\noindent and
$$\gamma_{k}(\sigma)(\cdots ,  t_p=0, \cdots) = $$ 
$$
\gamma_p(t_1,\cdots, t_{p-1})\gamma_{q}(t_{p+1},\cdots,t_{k-1}).$$

\begin{figure}[htbp]
\begin{center}
\includegraphics[width=.6\textwidth]{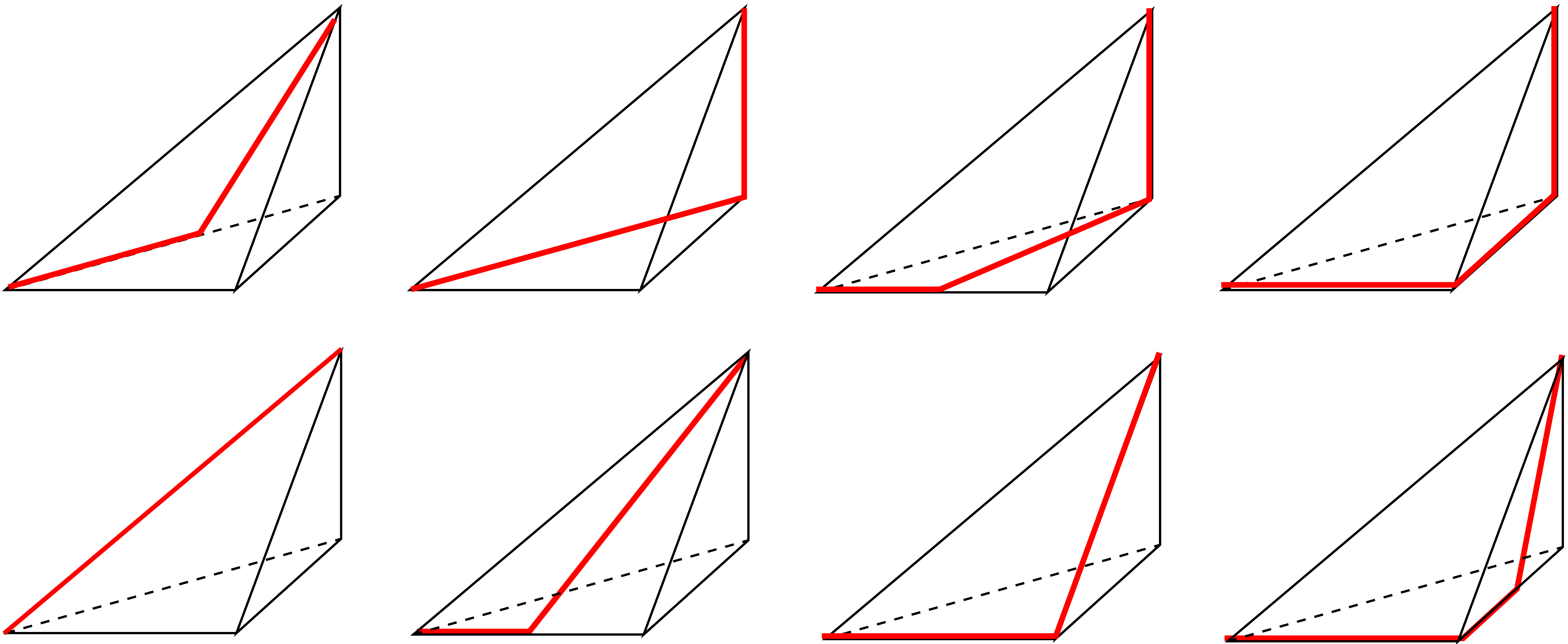}		
\caption{The relation between the $p_n$ and the $\gamma_n$}
\label{mine}
\end{center}
\end{figure}

One way to describe the relation between the $p_n$ and the $\gamma_n$  in words is: travel from vertex 0 partway toward vertex 1 then straight partway toward vertex 2 then 
straight partway toward vertex 3 etc. (See Figure\ref{mine}.)

Note that these are slighty different from the version of $\gamma_n$
given by Igusa; see  Figure \ref{igusa} taken from  \cite{igusa:iterated}.


\def\tetrahedron{{\thinlines\color{green}  
  \put(0,0){\line(1,0){.6}}  
  \qbezier(.6,0)(.7,.1)(.9,.3) 
  \put(.9,.3){\line(0,1){.6}} 
  \qbezier(0,0)(.3,.3)(.9,.9) 
  \qbezier(0,0)(.3,.1)(.9,.3)  
  \qbezier(.6,0)(.7,.3)(.9,.9)  
    }}
\def\tetrahedronB{{\thinlines\color{green}  
  \put(0,0){\line(1,0){.6}}  
  \qbezier(.6,0)(.7,.1)(.9,.3) 
  \qbezier(0,0)(.3,.3)(.9,.9) 
  \qbezier(0,0)(.3,.1)(.9,.3)  
  \qbezier(.6,0)(.7,.3)(.9,.9)  
    }}
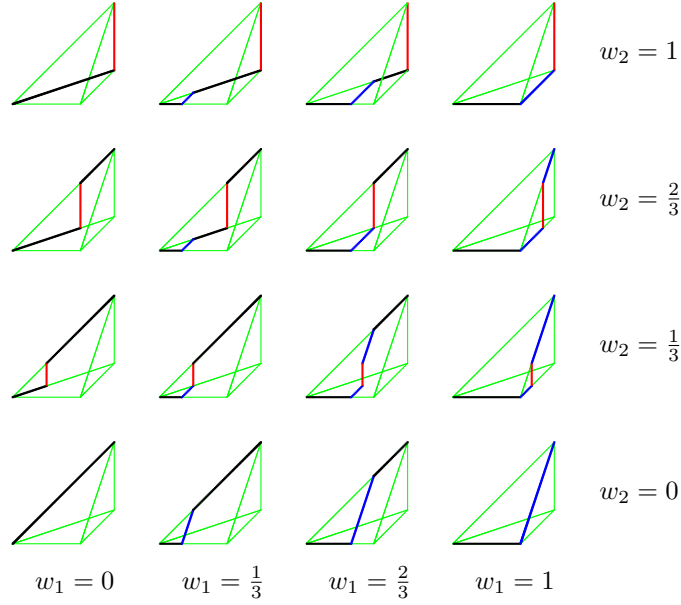
\begin{figure}[htbp]
\begin{center}
%
{
\setlength{\unitlength}{1.5cm}
{\mbox{
\begin{picture}(5.5,5.2)
	\put(0,.3){
		\thinlines
	\put(.2,-.4){$w_1=0$}
	\put(1.5,-.4){$w_1=\frac13$}
	\put(2.8,-.4){$w_1=\frac23$}
	\put(4.1,-.4){$w_1=1$}
	\put(5.2,.4){$w_2=0$}
	\put(5.2,1.7){$w_2=\frac13$}
	\put(5.2,3){$w_2=\frac23$}
	\put(5.2,4.3){$w_2=1$}
	%
    \put(0,0){\tetrahedron \thicklines\qbezier(0,0)(.3,.3)(.9,.9)}
    \put(0,1.3){\tetrahedron \thicklines\qbezier(.3,.3)(.7,.7)(.9,.9) 
     \qbezier(0,0)(.15,.05)(.3,.1) }
    \put(0,2.6){\tetrahedron \thicklines
     \qbezier(0,0)(.3,.1)(.6,.2)}
    \put(0,3.9){\tetrahedronB \thicklines
     \qbezier(0,0)(.3,.1)(.9,.3)}
    \put(1.3,0){\tetrahedron  \thicklines\put(0,0){\line(1,0){.2}} 
    {\color{blue}
     \qbezier(.2,0)(.25,.15)(.3,.3) } 
      \qbezier(.3,.3)(.7,.7)(.9,.9)  }
          \put(0,1.3){\thicklines {\color{red}\put(.3,.1){\line(0,1){.2}} }}
    \put(1.3,1.3){\tetrahedron  \thicklines\put(0,0){\line(1,0){.2}}
    {\color{blue}
     \qbezier(.2,0)(.25,.05)(.3,.1) }}
    \put(1.3,2.6){\tetrahedron  \thicklines\put(0,0){\line(1,0){.2}}
     \qbezier(.3,.1)(.45,.15)(.6,.2){\color{blue}
     \qbezier(.2,0)(.25,.05)(.3,.1) }}
    \put(1.3,3.9){\tetrahedronB  \thicklines\put(0,0){\line(1,0){.2}}
     \qbezier(.3,.1)(.6,.2)(.9,.3)
     {\color{blue}
     \qbezier(.2,0)(.25,.05)(.3,.1) }}
     \put(1.3,1.3){\thicklines {\color{red}\put(.3,.1){\line(0,1){.2}}}
      \qbezier(.3,.3)(.7,.7)(.9,.9) }
    \put(2.6,0){\tetrahedron \thicklines\put(0,0){\line(1,0){.4}} 
    {\color{blue}
     \qbezier(.4,0)(.5,.3)(.6,.6) }  
      \qbezier(.6,.6)(.7,.7)(.9,.9) }
    \put(2.6,1.3){\tetrahedron \thicklines\put(0,0){\line(1,0){.4}}
    {\color{blue}
     \qbezier(.4,0)(.45,.05)(.5,.1) \qbezier(.5,.3)(.55,.45)(.6,.6)}}
    \put(2.6,1.3){\thicklines {\color{red}\put(.5,.1){\line(0,1){.2}}
    } 
     \qbezier(.6,.6)(.7,.7)(.9,.9)}
    \put(2.6,2.6){\tetrahedron \thicklines\put(0,0){\line(1,0){.4}}
    {\color{blue}
     \qbezier(.4,0)(.5,.1)(.6,.2) }  }
    \put(0,2.6){\thicklines {\color{red}\put(.6,.2){\line(0,1){.4}}
    }
    \qbezier(.6,.6)(.7,.7)(.9,.9) }
    \put(1.3,2.6){\thicklines {\color{red}\put(.6,.2){\line(0,1){.4}}} 
     \qbezier(.6,.6)(.7,.7)(.9,.9)}
    \put(2.6,2.6){\thicklines {\color{red}\put(.6,.2){\line(0,1){.4}}}
     \qbezier(.6,.6)(.7,.7)(.9,.9) }
    \put(0,3.9){\thicklines {\color{red}\put(.9,.3){\line(0,1){.6}}
    } }
   \put(1.3,3.9){\thicklines {\color{red}\put(.9,.3){\line(0,1){.6}}} }
   \put(2.6,3.9){\thicklines {\color{red}\put(.9,.3){\line(0,1){.6}}} 
   }
   \put(2.6,3.9){\thicklines
    \qbezier(.6,.2)(.75,.25)(.9,.3)}
   \put(2.6,3.9){{\thinlines\color{green}  
  \put(0,0){\line(1,0){.6}}  
  \qbezier(.6,0)(.7,.1)(.9,.3) 
  \qbezier(0,0)(.3,.3)(.9,.9) 
  \qbezier(0,0)(.3,.1)(.6,.2)  
  \qbezier(.6,0)(.7,.3)(.9,.9)  
    } \thicklines\put(0,0){\line(1,0){.4}}
     {\color{blue}
     \qbezier(.4,0)(.5,.1)(.6,.2) } }
     \put(3.9,0){\tetrahedron \thicklines\put(0,0){\line(1,0){.6}}
     {\color{blue}
     \qbezier(.6,0)(.7,.3)(.9,.9) } }
    \put(3.9,1.3){\tetrahedron \thicklines\put(0,0){\line(1,0){.6}}
    {\color{blue}
     \qbezier(.6,0)(.65,.05)(.7,.1)\qbezier(.7,.3)(.8,.6)(.9,.9) }}
    \put(3.9,1.3){\thicklines {\color{red}\put(.7,.1){\line(0,1){.2}} }}
     \put(3.9,2.6){\tetrahedron \thicklines\put(0,0){\line(1,0){.6}} 
     {\color{blue}
     \qbezier(.6,0)(.7,.1)(.8,.2)\qbezier(.8,.6)(.85,.75)(.9,.9) }}
     \put(3.9,2.6){\thicklines {\color{red}\put(.8,.2){\line(0,1){.4}} 
     } }
   \put(3.9,3.9){\tetrahedronB}
    \put(3.9,3.9){
   \thicklines\put(0,0){\line(1,0){.6}}
    {\color{blue}\qbezier(.6,0)(.7,.1)(.9,.3)}}
	 \put(3.9,3.9){\thicklines{\color{red}\put(.9,.3){\line(0,1){.6}}
	 }  }}
\end{picture}
}}
\caption{Igusa's Figure 3}
\label{igusa}
}
\end{center}
\end{figure}
%


Hopefully the pattern is clear.

Correspondingly, the liftings $\Theta_n:I^n \times E\to E$ form a 
collection of maps which assign to any $k$-simplex $\sigma:\Delta^k\to B$ a map $ \Theta_{k}(\sigma): I^{k}\times F_0 \to F_k$ satisfying, for any $e\in F_0$, the relations :

$\Theta_0(0)$ is the identity on  $F_0.$

\noindent For any $(t_1,\cdots,t_{k})$, 

$\Theta_k(\sigma)(t_1,\cdots,t_{k}, -): F_0 \to F_k$ is a homotopy equivalence.

\noindent For any $1\leq p\leq k-1$,

$$ \Theta_{k}(\sigma)(\cdots ,  t_p=0, \cdots, e) =  \Theta_{k-1}(\partial_p\sigma)(\cdots,\hat t_p,\cdots, e)$$

$$ \Theta_{k}(\sigma)(\cdots ,  t_p=1, \cdots, e) = $$ 
$$
 \Theta_{p}(<0,\cdots,p>)(t_1,\cdots, t_{p-1},  \theta_{q}(<p,\cdots,k>)(t_{p+1},\cdots,t_k,e)).$$
\vskip1ex
\noindent The desired $\theta_n$ is again recovered at $t_1=1.$
\vskip2ex
The maps $p_n$ can be interpreted as homotopies $q_n: I\to (\Delta^n)^{I^{n-1}}$ and so subject to the homotopy lifting property.
For example, $\gamma_1: 0\to P\Delta^1$ is a path which can be lifted as in Theorem 1 to give 
$\Theta_1: I\times F_0\to E$. Then
$\gamma_2:I\to P\Delta^2$ such that $0$ maps to the `identity' path $I\to <02>$ while $1$ maps to the concatenated path $<01><12>.$ (Henceforth, we will assume paths have been normalized to length 1 where appropriate.) Now lift the homotopy $\gamma_2$ to a homotopy 
$\Theta_2(<012>):I\times I\times E \to E$ between $\Theta_1( <02>)$ and $\Theta_1(<01><12>).$ In particular,
$\Theta_2(<012>):1\times I\times E \to E$ gives the desired homotopy $\theta_2:I\times F_0\to F_k.$

The situation becomes slightly more complicated as we increase the dimension.
The case 
$\Delta^3$ is illustrative. The faces $<023>$ and $<013>$ lift just as $<012>$ had via $\Theta_2$, but 
that lift must then be `whiskered'  by a rectangle over $<23>$ which glues onto $\Theta_3(<012>$. In a less complicated way $<123>$ is lifted so that
vertex $1$ agrees with the end of the `whisker' which is the lift of $<01>$. Thus the total lift of $<0123>$ ends with the desired $\theta_3: I^2\times F_0\to F_3.$ The needed whiskering (of various dimensions) is prescribed by the $t_p=1$ relations of Definition \ref{Thetas} to be satisfied. 

\begin{figure}[htbp]
\begin{center}
\includegraphics[width=.6\textwidth]{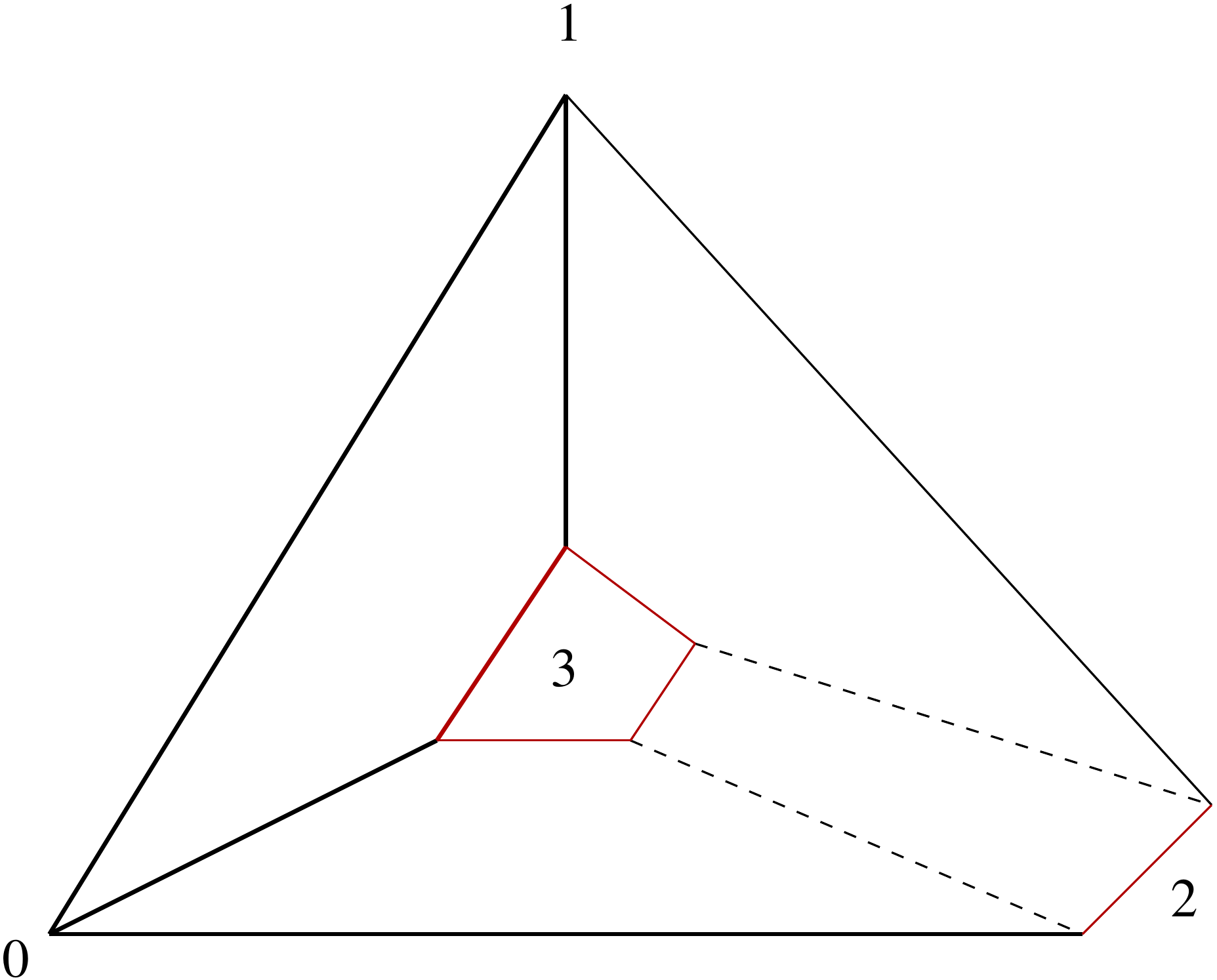}		
\caption{$\Theta_3$} 
\label{Theta3}
\end{center}
\end{figure}
 Further details are left to the industrious reader.
\section{(Re)-construction of fibrations}
In \cite{jds:mex}, I showed how to construct a fibration from the data of an strong homotopy action  of  
$\Omega B$ on a `fibre' $F$. If the action came from a given fibration $F\to E\to B$, the constructed fibration
was fibre homotopy equivalent to the given one. For \emph{representations up to homotopy}, a similar result applies using analogous techniques, with some additional subtlety.

First we try to construct a fibration naively. Over each 1-simplex $\sigma$ of $Sing(B)$, we take 
$\sigma\times F_0$ and attempt to glue these pieces appropriately. For the one simplices $<01>$ and $<12>$,
we have $\theta_1:F_0\to F_1$ which tells us how to glue $<01>\times F_0$ to $<12>\times F_1$ at vertex 1,
but, since $\theta_1:F_0\to F_2$ is not the composite of  $\theta_1:F_0\to F_1$ and  $\theta_?:F_1\to F_2$,
we can not simply plug in $<012>\times F_0$ over  $<012>$.  However, we can plug in $I^2\times F_0$
since $\theta_2: I\times F_0\to F_2$ will supply the glue over vertex 2.

To describe the total space of the fibration (or at least a quasi-fibration),  we use the special maps $p_n:I^n\to \Delta^n.$ 
Now return to the description of the fibration $\bar p_2: E_2\to\Delta^2$ above.
  In greater precision,
  $$E_2= <01>\times F_0\ \cup_1 <12>\times F_1\ \cup_0 <02>\times F_0\  \cup I^2\times F_0.$$
  The attaching maps over the vertices $0$ and $1$ are obvious as are the projections to the edges of 
  $\Delta^2.$  On $ I^2\times F_0$, the attaching maps are obvious except for the face $t_2=s=1$ where it is given by $\theta_2: I\times F_0\to F_2,$
  so as to be compatible with the projection $\bar p_2: I^2\times F_0 \to \Delta^2.$
  
    The result is at least a quasi-fibration $q:E_\theta\to B$ and can be replaced up to fibre homotopy equivalence by a true fibration (cf. \cite{ fuchs:fibrations}.
  
 Notice that although the definition of representation up to homotopy was in terms of  fibrations, in fact it really needs only the collection of fibres $F_\sigma$ for the 0-simplices of $Sing(B)$. The equivalence in the appropriate sense between \emph{representations up to homotopy of $Sing(B)$} and \emph{fibrations over $B$} follows as for Theorem B in \cite{jds:mex}.
 
 \section{Coda}
 Physics is often written in terms of smooth structures, differential forms and `geometrically' in terms of connections.
 From a topological point of view, parallel transport is the more basic notion. In particular, string theory and string field theory
 has inspired string topology, inititiated by Chas and Sullivan,  and a variety of $\infty$-algebras. I look forward to the corresponding representation-up-to-strong homotopy theory feeding back into physics.


\begin{thebibliography}{10}

\bibitem{abad-schaetz}
C.~Abad and F~Schatz.
\newblock The $\mathrm{A}_\infty$ de {R}ham theorem and integration of
  representations up to homotopy.
\newblock in progress.

\bibitem{adams:cobar}
J.~F. Adams.
\newblock On the cobar construction.
\newblock {\em Proc. Nat. Acad. Sci. USA}, 42:409--412, 1956.

\bibitem{block-smith}
J.~Block and A.~Smith.
\newblock A {R}iemann {H}ilbert correspondence for infinity local systems.
\newblock {\tt{ arXiv:0908.2843}}.

\bibitem{chen:iterated}
K.-T. Chen.
\newblock Iterated integrals of differential forms and loop space homology.
\newblock {\em Annals of {M}athematics}, 97:217--246, 1973.

\bibitem{duran:exotic}
C.~E. Dur{\'a}n.
\newblock Pointed {W}iedersehen metrics on exotic spheres and diffeomorphisms
  of {$S^6$}.
\newblock {\em Geom. Dedicata}, 88(1-3):199--210, 2001.

\bibitem{d&e}
Eldon Dyer and Samuel Eilenberg.
\newblock Globalizing fibrations by schedules.
\newblock {\em Fund. Math.}, 130(2):125--136, 1988.

\bibitem{fss}
D.~Fiorenza, U.~Schreiber, and J.~Stasheff.
\newblock Cech cocycles for differential characteristic classes -- an
  infinity-lie theoretic construction.
\newblock {\tt{arXiv:1011.4735}}, 2010.

\bibitem{fuchs:fibrations}
M.~Fuchs.
\newblock A modified {D}old-{L}ashof construction that does classify
  {$H$}-principal fibrations.
\newblock {\em Math. Ann.}, 192:328--340, 1971.

\bibitem{gugenheim:deRham}
V.~K. A.~M. Gugenheim.
\newblock On the multiplicative structure of the de {R}ham theory.
\newblock {\em J. Differential Geometry}, 11(2):309--314, 1976.

\bibitem{hilton:intro}
P.~J. Hilton.
\newblock {\em An introduction to homotopy theory}.
\newblock Cambridge Tracts in Mathematics and Mathematical Physics, no. 43.
  Cambridge, at the University Press, 1953.

\bibitem{igusa:iterated}
K.~Igusa.
\newblock Iterated integrals of superconnections.
\newblock {\tt{arXiv:0912.0249}}.

\bibitem{igusa:twisting}
K.~Igusa.
\newblock Twisting cochains and higher torsion.
\newblock {\tt{arXiv:math/0212383 }}.

\bibitem{iwase-mimura:hha}
Norio Iwase and Mamoru Mimura.
\newblock Higher homotopy associativity.
\newblock In {\em Algebraic topology ({A}rcata, {CA}, 1986)}, volume 1370 of
  {\em Lecture Notes in Math.}, pages 193--220. Springer, Berlin, 1989.

\bibitem{MS:cc}
J.~Milnor and J.~Stasheff.
\newblock {\em Characteristic {C}lasses}, volume~76 of {\em Annals of Math
  Studies}.
\newblock Princeton Univ. Press, 1974.

\bibitem{jcm:paths}
J.~C. Moore.
\newblock Le th\'eor\`eme de {F}reudenthal, la suite exacte de {J}ames et
  l'invariant de {H}opf g\'en\'eralis\'e.
\newblock In {\em S\'eminaire Henri Cartan;1954-55}, pages 22--01 -- 22--??
  CBRM, 1955.

\bibitem{nowlan:action}
R.~A. Nowlan.
\newblock {$A_{n}$}-actions on fibre spaces.
\newblock {\em Indiana Univ. Math. J.}, 21:285--313, 1971/1972.

\bibitem{jds:hahI}
J.~Stasheff.
\newblock Homotopy associativity of {H}-spaces, {I}.
\newblock {\em Trans. Amer. Math. Soc.}, 108:293--312, 1963.

\bibitem{jds:shrep}
J.~Stasheff.
\newblock Constrained {P}oisson algebras and strong homotopy representations.
\newblock {\em Bull. Amer. Math. Soc. (N.S.)}, 19(1):287--290, 1988.

\bibitem{jds:madison}
James~D. Stasheff.
\newblock {$H$}-spaces and classifying spaces: foundations and recent
  developments.
\newblock In {\em Algebraic topology ({P}roc. {S}ympos. {P}ure {M}ath., {V}ol.
  {XXII}, {U}niv. {W}isconsin, {M}adison, {W}is., 1970)}, pages 247--272. Amer.
  Math. Soc., Providence, R.I., 1971.

\bibitem{jds:bingh}
James~D. Stasheff.
\newblock Parallel transport and classification of fibrations.
\newblock In {\em Algebraic and geometrical methods in topology (Conf.
  Topological Methods in Algebraic Topology, State Univ. New York, Binghamton,
  N.Y., 1973)}, pages 1--17. Lecture Notes in Math., Vol. 428. Springer,
  Berlin, 1974.

\bibitem{jds:mex}
James~Dillon Stasheff.
\newblock ``{P}arallel'' transport in fibre spaces.
\newblock {\em Bol. Soc. Mat. Mexicana (2)}, 11:68--84, 1966.
\newblock {\tt{scanned copy available from author}}.

\bibitem{sugawara:g}
M~Sugawara.
\newblock On a condition that a space is group-like.
\newblock {\em Math. J. Okayama Univ.}, 7:123--149, 1957.

\bibitem{sugawara:hc}
M~Sugawara.
\newblock On the homotopy-commutativity of groups and loop spaces.
\newblock {\em Mem. Coll. Sci. Univ. Kyoto, Ser. A Math.}, 33:257--269,
  1960/61.

\end{thebibliography}
\end{document}